\documentclass[a4paper,10pt]{article}
\usepackage{epsfig}
\usepackage{amsmath}
\usepackage{amssymb}
\usepackage{amscd}

\begin{document}
\title{ON THE CLASSIFICATION OF ONE CLASS OF THREE-DIMENSIONAL BOL ALGEBRAS WITH SOLVABLE ENVELOPING LIE ALGEBRAS OF SMALL DIMENSION I}
\author{THOMAS BOUETOU BOUETOU}

\maketitle \begin{quote}{\bf Abstract.}
\em We give the classification of 3-dimensional Bol algebras associated to the
trilinears operations  $(e_2,e_3,e_2)=e_1 $ ; \;  $(e_2,e_3,e_3)=\pm e_2 $. It
 turn that for such trilinears operation there exist up to isomorphism and
 isotopy one family and one exceptional Bol algebras. We describe their associated 3-webs.
\end{quote}

\maketitle \begin{quote}{\bf Keywords:}
\em Bol algebras, Lie groups, Lie algebras, Lie triple system, B-Isotopic Bol algebras, analytic Bol loops, Bol 3-Webs, 
\end{quote}
\maketitle \begin{quote}{\bf AMS subject classification 2000:}
 11E57, 14C21, 17D99, 19L99, 20N05, 22E60, 22E67, 32C22, 53A60
\end{quote}

\begin{center}

\section{INTRODUCTION}

\end{center}

The notion of quasigroup and the tangential structure appearing to it
have been intensively developing for the 40 past years. The quasigroup as
a nonassociative algebraic structure is naturally the generalization
of the notion of group. It first appeared in the work of R. Moufang
(1935) \cite{mou}. She obtained some identities. The smooth local loops were
first introduced in the work of Malcev A.I. \cite{malc}, in connection with
the generalization of Lie group.
 
Left binary algebras were later on called Malcev algebras. It is well
known that modern differential geometry and nonassociative algebras
are interacting on one another. The notion of binary-ternary
operation tangent to the given geodesic loop connected with an
arbitrary point in the affine connected space was introduced in Akivis
works see \cite{ak1,ak2}. Several mathematicians worked in the development of
differential geometry and the study of different classes of
quasigroups and loops \cite{loos,sabmik4,fed2,shele1,shele2,bar,bru,kik1}.

In 1925 E. Cartan in his research lay the beginning of investigation
of symmetric spaces. Today the given structure plays an important role,
in differential geometry and its application. The question arising
here is that of, the description and the classification of symmetric
space naturally leading to the classification of the corresponding
algebraic structure. 

The survey on geometry of fiber space stimulates the interest for
special types of 3-webs, in particular Bol 3-Webs and the tangential
structures appearing to it: Bol algebras \cite{mike2,mike6}. In 
connection with it
the idea of the description of collections of three-dimensional Bol
algebras elaborates. For separate classes of Bol algebras see 
\cite{alhs,fedo}. In
this investigation our approach on Bol algebras will be base on the
classification of solvable Lie triple systems \cite{boue1}. As we know Bol
algebra can be seen as a Lie triple system equipped with an additional
bilinear anti-commutative operation which verify a series of
supplementary conditions. 

Bol algebras appear under the infinitesimal description of the class
of local smooth Bol loop, in the work of L.V. Sabinin and
P.O. Mikheev. The interest of studying smooth Bol loop is connected to
the fact that, the geodesical loop see(Sabinin \cite{sab4}) of local
 symmetric
affine connected space, verifies the left Bol identity and automorphic
inverse identity. Its an exact algebraic analog of construction of
symmetric spaces. In particular, the velocity space in the theory of
relativity (STR), is a Bol loop relatively to the additional law of
velocity \cite{sabmik5}. Thus, relying on what is stated above we 
are ready to
formulate the purpose of this work: Classification of Bol algebras of
dimension 3 with solvable Lie algebras of dimension $\leq
5$, with accuracy to isomorphism and with accuracy to isotopic which
include: 

\begin{enumerate}
\item the  description of three-dimensional Lie triple systems and
their corresponding Lie algebras with invomorphisme. 
\item The description of three-dimensional Bol algebras linked with
the distinguished Lie triple systems above. 
\item The construction of the object describing the isotopy of Bol algebras.
\item The description of Bol 3-Webs, connected with the selected Bol algebras.
\end{enumerate} 

\section{Basic definition}

{\bf Definition 1.} A Bol algebra is called a vector space $V$ over the field of characteristic $0$  with the two operations 
$$\xi,\eta \to \xi\cdot\eta\in V,$$
$$\xi,\eta,\zeta\to (\xi,\eta,\zeta)\in V \quad (\xi,\eta,\zeta\in
V)$$ satisfying the following  identities 
\begin{enumerate}
\item $$ \xi\cdot\xi=0,\quad
(\xi,\eta,\zeta)+(\eta,\zeta,\xi)+(\zeta,\xi,\eta)=0,$$

\item $$(\xi,\eta,\zeta)
\chi-(\xi,\eta,\chi)\zeta+(\zeta,\chi,\xi\cdot\eta)-
(\xi,\eta,\zeta\cdot\chi)+\xi\eta\cdot\zeta\chi=0,$$
\item $$(\xi,\eta,(\zeta,\chi,\omega))=((\xi,\eta,\zeta),\chi,\omega)+(\zeta
,(\xi,\eta,\chi),\omega)+(\zeta,\chi,(\xi,\eta,\omega)).$$

\end{enumerate}

Under the base of the  definition above the following lemma holds:

{\bf Lemma 1}  Every vector space $ \mathfrak{B} $   whose main operations are
 $\kappa \cdot \zeta $ and
$(\xi, \eta, \chi) $ such that there exist a linear mapping 
$$\prod : \mathfrak{B} \wedge  \mathfrak{B}  \longrightarrow \mathcal{End}(\mathfrak{B})$$
$$(\xi, \eta)\longmapsto \prod _{\xi \eta}:  \mathfrak{B} \longrightarrow \mathfrak{B}$$
 with companion 
$\xi \cdot \eta $ define by $\prod _{\xi,\eta}(\chi)=(\xi, \eta, \chi)$
which is a pseudo-derivative on the bilinear operation and a derivative on a trilinear operation. that is:
\begin{itemize}
\item $$\prod (\kappa \cdot \zeta)=\prod \kappa \cdot \zeta+\kappa \cdot \prod \zeta+(\kappa, \zeta, \xi \cdot \eta)-\kappa \zeta \cdot \xi \eta$$

\item $$\prod (\xi, \eta, \chi)=(\prod \xi, \eta, \chi,)+(\xi, \prod \eta, \chi)+(\xi, \eta, \prod \chi)$$
\end{itemize}
for any $\xi, \eta, \zeta, \chi, \kappa $, is Bol algebra that is verify condition 1-3 of definition 1.

Since a Bol algebra can be seen as a Lie triple system equipped with additional
identities, this fact can be caracterise by the following lemma:

{\bf Lemma 2.} A Lie triple system $\mathfrak{M}$ with a skew symmetric bilinear 
operator $(\cdot)$ is a Bol algebra if $\forall \xi, \eta \in \mathfrak{M}$
the mapping $\prod $ defined by:
$$\prod_{\xi, \eta} : \mathfrak{M} \longrightarrow \mathfrak{M}$$ 
$$\chi \longmapsto (\xi, \eta, \chi)$$

is a pseudo-derivation of the bilinear operation 
$(\cdot)$ with $\xi \cdot \eta$ as companion and a derivation of the trilinear
operation.

That is:

 $\forall \kappa, \zeta \in \mathfrak{M}$

$$\prod_{\xi, \eta} (\kappa \cdot \zeta)=(\prod_{\xi, \eta} \kappa) \cdot \zeta+\kappa \cdot (\prod_{\xi, \eta} \zeta)+(\kappa, \zeta, \xi \cdot \eta)-\kappa \zeta \cdot \xi \eta$$

$$\prod (\xi, \eta, \chi)=(\prod \xi, \eta, \chi,)+(\xi, \prod \eta, \chi)+(\xi, \eta, \prod \chi).$$

The proof of this lemma is based under a direct computation.

We will note that any Bol algebra can be realized as the tangent 
algebra to a  Bol Loop with the left Bol identity, and they allow embedding in Lie
algebras. 

This can be expressed this way;
let $ (G, \Delta, e) $- be a local Lie group,  $ H $-  one of its subgroups, 
and let us  denote the  corresponding Lie algebra and subalgebra by
$ \mathfrak{G}$ and $ \mathfrak{h}$. 

Consider a vector subspace 
$\mathfrak{B}$ such that 

$ \mathfrak{G}=\mathfrak{h}\dotplus \mathfrak{B}$.
 Let $ \Pi:G\longrightarrow G\setminus H $ be the canonical projection and let
$ \Psi$ be the restriction of mappings composition $ \Pi \circ \exp$, 
to $ \mathfrak{B}$. Then there exists such a neighborhood $ \mathcal{U}$ of 
the
point $ O $ in $\mathfrak{B}$ such that $\Psi$ maps it diffeomorphically into the 
neighborhood $\Psi (u)$ of the coset $\Pi (e)$ in $ G \setminus H $ .\\

$$
\begin{CD}
\mathfrak{G} @<{i}<< \mathfrak{B}\\
@V{\exp}VV  @VV{\Psi}V \\
G @>{\pi}>> G\setminus H
\end{CD}
$$

by introducing a local composition law:

$$
 a \star b= \Pi_{B}(a\Delta b) 
$$ 
 
On points of local cross-section $ B=\exp \mathcal{U} $ of left coset of $G$
 mod $H$ where $ \prod_{B}=\exp \circ \Psi^{-1} \circ\Pi: G \longrightarrow B $
is the local projection on $ B $ parallel to $ H $ which puts  every element 
$ a \in B $ in correspondence so that $ g=a \Delta p $, where $  p \in H $.
We will emphasize that if any two local analytic Loops are isomorphic then 
their
corresponding Bol algebras are isomorphic.

\underline{{\bf Proposition 1.}}[Sabinin-Mikheev] Let us assume that for any 
$ a, b \in B=\exp \mathcal{U} $ sufficiently close to the point $e$, and 
$ a\Delta b \Delta a \in B $; then the local analytic loop 
$ (B, \times, e) $ satisfies the left Bol condition.


\section{EMBEDDING OF BOL ALGEBRA INTO LIE ALGEBRAS: ENVELOPING LIE ALGEBRA OF BOL ALGEBRAS}


Let us suppose that the local cross-section $ B=\exp \mathcal{U} $ of left cosets $ G
\setminus H $ satisfy the condition of Proposition 1. above. It is
interesting to calculate the operations of the Bol algebra $
\mathfrak{B} $ tangent to the local analytic Bol loop $ ( B, \star, e)$
in terms of a Lie algebra $ \mathfrak{G} $, its subalgebra 
$\mathfrak{h} $ and vector subspace $ \mathfrak{B} $. Let us introduce in $ G $
normal coordinates then: 

$$
a \star b=a+b+\frac{1}{2}[a,b]_{\mathfrak{B}}+0(2)
$$
$\forall a $ ,$ b \in B $
where $ [a,b]_{\mathfrak{B}} $ is the projection of $ [a,b] $ on $ \mathfrak{B} $ parallel to $ \mathfrak{h} $ thus
$$
\xi \cdot \eta=[\xi, \eta]_{\mathfrak{B}}
$$
$$
(\xi, \eta, \chi)=[[\xi, \eta], \chi]
$$
$$
<\xi, \eta, \chi>=-\frac{1}{2}[[\xi, \eta], \chi]+\frac{1}{2}[[\xi, \eta]_{\mathfrak{B}}, \chi]_{\mathfrak{B}}.
$$
One can find the correctness of the following propositions [30,44]:\\

\underline{{\bf Proposition 2.}} \cite{sabmik1} Two local analytic Bol 
loops  
$ B $ and 
$ B' $ are locally isomorphic, if and only if their corresponding Bol algebras 
$ \mathfrak{B} $ and $ \mathfrak{B'} $ are isomorphic.

\underline{{\bf Proposition 3.}} \cite{sabmik1} Every Bol algebra is
 a tangent 
$ W $-algebra of some local analytic loop with the left Bol loop.

\underline{{\bf Proposition 4.}} \cite{sabmik1}  The property 
$ a\Delta b \Delta a \in B $ 
is held for
 any $ a$, $b \in B $ sufficiently close to $e $, if and only if
 $ [[\xi, \eta], \zeta ] \in \mathfrak{B} 
 \forall \xi, \eta, \zeta \in \mathfrak{B} $.

Given a finite dimensional Bol algebra $\mathfrak{B}$  over $ \mathbb{R}$, 
the basic
 operations of which are $\xi \cdot \eta $ and $ (\xi, \eta, \zeta) $, we 
plane showing that there exist a finite dimensional Lie algebra
$ \mathfrak{G} $  over $\mathbb{R} $, containing a Lie subalgebra 
$\mathfrak{h}$   and a linear mapping 
$i$:$ \mathfrak{B} \longrightarrow \mathfrak{G} $ verifying 
 $ i(\mathfrak{B}) \in \mathfrak{B}$ such that
 .

$ \mathfrak{G}= \mathfrak{B}\dotplus \mathfrak{h} $ (direct sum of vector spaces)\\

$ [[\mathfrak{B},\mathfrak{B}], \mathfrak{B}] \in \mathfrak{B} $ and
$\forall \xi, \eta, \zeta \in \mathfrak{B}$ 

$ \xi \cdot \eta=[\xi, \eta]_{\mathfrak{B}} $
$ (\xi, \eta, \zeta)=[[\xi, \eta], \zeta] $

where $[\xi, \eta]$ denotes the result of commutation of vectors in $ \mathfrak{G}$
 and $ [ \xi, \eta]_{\mathfrak{B}} $ denotes projection of the vector
 $ [\xi, \eta] $ on $\mathfrak{B}$ parallel to $\mathfrak{h}$.

 In that case we will talk about the enveloping pair $ (\mathfrak{G}, \mathfrak{h})$ Lie algebra of Bol algebra$ \mathfrak{B}$ or, in other words enveloping
 Lie algebra $\mathfrak{G}$ of Bol algebra $\mathfrak{B}$.

Let   $ (\mathfrak{G}, \mathfrak{h}) $ be an enveloping pair of Lie algebra of
 Bol algebra$ \mathfrak{B}$. Let us identify $\mathfrak{B}$ with a vector
 subspace
 $ i(\mathfrak{B})$ into $\mathfrak{G}$, and let us consider the subalgebra
 $ \mathfrak{G'}=\mathfrak{B}\dotplus [ \mathfrak{B}, \mathfrak{B}] $ into
 $\mathfrak{G}$, and the subalgebra
 $ \mathfrak{h'}=\mathfrak{h}\cap \mathfrak{G'} $ into$ \mathfrak{G'}$, then 
the pair $ ( \mathfrak{G'}, \mathfrak{h'} ) $ is also for enveloping for a
 Bol algebra $\mathfrak{B}$.

By the construction of the Lie algebra which is a canonical enveloping for
$ \mathfrak{B} $, it is better to use the construction made in \cite{sbl}.
 For such
 Bol algebra $ \mathfrak{B}$, there exists a Lie algebra
 $  \mathfrak{\widetilde{G}} $, and an envoluting automorphism
 $ \tau \in Aut \mathfrak{\widetilde{G}} $ such that $ ( \tau^{2}=Id) $, linear
 injection map $ i: \mathfrak{B} \longrightarrow \mathfrak{\widetilde{G}} $
 and a subalgebra $ \mathfrak{\widetilde{h}} $ in $ \mathfrak{\widetilde{G}} $,
 such that ( we are identifying $ i(\mathfrak{B}) $ with $\mathfrak{B}$).

$ \mathfrak{\widetilde{G}} =\mathfrak{\widetilde{G_{-}}}+\mathfrak{\widetilde{G_{+}}} $, where $ \mathfrak{\widetilde{G_{-}}}=\mathfrak{B} $;
 $ \mathfrak{\widetilde{G}}=\mathfrak{B} \dotplus \mathfrak{\widetilde{h}} $,
 $ < \mathfrak{B}>=\mathfrak{\widetilde{G}} $

$$
\xi \cdot \eta= \prod_{\mathfrak{B}}[\xi, \eta]= [\xi, \eta]_{\mathfrak{B}}
$$
$$
(\xi, \eta, \zeta)=[[\xi, \eta], \zeta]
$$

Lie algebra $ \mathfrak{\widetilde{G}}$ is an enveloping Lie algebra for Bol 
algebra $ \mathfrak{B}$, but in general not canonical enveloping because 
$ \mathfrak{\widetilde{h}}$ may contain an ideal $ I $ of Lie algebra
$ \mathfrak{\widetilde{G}}$. Hence the canonical enveloping Lie algebra 
$ \mathfrak{G}$ for a Bol algebra $ \mathfrak{B}$ is obtained by factorizing
$ \mathfrak{\widetilde{G}}$  with the ideal $ I $. Therefore after the 
factorization of $ \mathfrak{\widetilde{G}}_{+} \setminus I $ and $\mathfrak{B}$
 ( we identify $\mathfrak{B}$ and $ \mathfrak{B} \setminus I $) in general 
interest, let us note that the construction see \cite{sbl} follows that 
$ dim \mathfrak{\widetilde{G}} \leq dim \mathfrak{B} \wedge \mathfrak{B}+dim \mathfrak{B} $.

In our case $ dim \mathfrak{B}=3 $, that is why under the examination of the
 corresponding canonical enveloping Lie algebra $\mathfrak{G}$ we must consider
 the case:
$$
dim[ \mathfrak{B}, \mathfrak{B}]=0, 1, 2, 3.
$$  

 This embedding can be justify also this way:

Let $(\mathfrak{B}, (\cdot ), (, , ))$ be a Bol algebra. Embedding $\mathfrak{B}$
in a Lie algebra  $(\mathfrak{M}$ mean we can find  $(\mathfrak{h}$ such that
 $\mathfrak{G}= \mathfrak{B}\bigoplus  \mathfrak{h}$ and  $(\mathfrak{h}$ is 
a subalgebra of  $(\mathfrak{G}$.
we have 
$$ \mathfrak{h} \longrightarrow  \mathfrak{G}\longrightarrow \mathfrak{G}/ \mathfrak{h} \backsimeq \mathfrak{B}$$

$$0 \longrightarrow \mathfrak{h} \longrightarrow  \mathfrak{G}\longrightarrow \mathfrak{G}/ \mathfrak{h} \backsimeq \mathfrak{B}  \longrightarrow 0$$

which is an exact sequences of $\mathfrak{h}-module$.

\section{ABOUT B-ISOTOPIC BOL ALGEBRAS}

Below we give a  generalization of the notion of isotopy of (global) 
loops to the case of local analytic Bol loops. Our approach is based on a 
construction of embedding of Bol loops into a group and on an interpretation 
of isotopic loops in terms of their enveloping groups \cite{sabmik3,mike8}.\\

\underline{{\bf Definition  2.}} Let $ \mathfrak{B}$ and $\mathfrak{\widetilde{B}} $ be 
two Bol algebras and let $ \mathfrak{G}=\mathfrak{B} \dotplus \mathfrak{h} $ 
and $ \mathfrak{\widetilde{G'}}=\mathfrak{\widetilde{B}}\dotplus \mathfrak{\widetilde{h}}$ 
be their canonical enveloping Lie algebras. The algebras $\mathfrak{B}$ and 
$ \mathfrak{\widetilde{B}}$ will be called isotopic if there exit such a Lie 
algebras isomorphism $ \Phi: \mathfrak{G} \longrightarrow \mathfrak{\widetilde{G}} $, 
such that $ \Phi(\mathfrak{G})=\mathfrak{\widetilde{G}} $ and $ \Phi (\mathfrak{h}) $ 
coincide with the image of the subalgebra $\mathfrak{\widetilde{h}}$ in
$ \mathfrak{\widetilde{G}}$ under the action of an inner automorphism
 $ Ad \xi $, $ \xi \in \mathfrak{\widetilde{G}}$ i.e.
$$
 \Phi (\mathfrak{h})=( Ad_{\xi}) \mathfrak{\widetilde{h}}
$$

It is clear that the notion of isotopy is not an equivalence relation to Bol algebra manifolds.\\

{\bf \underline{\bf Theorem 1.} (MIKHEEV-BOUETOU)} Let $ B(\times) $ and $ \widetilde{B}(\circ) $ 
be global analytic Bol loops, and let their tangent Bol algebras be isotopic, then $ \widetilde{B}(\circ) $ is locally isomorphic to an analytic Bol loop analytically isotopic to $ B(\times) $.\\

Proof.\\
Let  $ \mathfrak{G} $ be the canonical Lie algebra enveloping the Bol algebra  $ \mathfrak{B} $ tangent to the analytic Bol loop $ B(\times) $.  $ \mathfrak{G}=\mathfrak{B} \dotplus \mathfrak{h} $ ( direct sum of vector spaces). There exists a Lie group $ G $ with $ \mathfrak{G} $ as its tangent Lie algebra, a closed subgroup $ G_{O} $ corresponding to the subalgebra $ \mathfrak{h} $ and an analytical embedding $ i: B \longrightarrow G $ such that the composition law takes the following form \cite{sabmik3,mike8}
$$
a \times b=\prod(a,b),
$$
where $a \times b$ denotes the composition of elements $a$, $ b $ in $ B $ and $ \prod: G \longrightarrow B $ is the projection on $ B $ parallel to  $ G_{O} $.\\ 

Let us suppose that the Bol algebra $ \widetilde{\mathfrak{B}} $ which is tangent to the loop $ \widetilde{B}(\circ) $ is isotopic to the Bol algebra $ \mathfrak{B} $, i.e. there exists a  Lie algebra isomorphism
 $ \Phi : \mathfrak{G} \longrightarrow  \mathfrak{\widetilde{G}}$, such that  $ \Phi ( \mathfrak{G})=  \mathfrak{\widetilde{G}} $ and  $ \Phi ( \mathfrak{h}) $ coincides with the image of the subalgebra $  \mathfrak{\widetilde{h}}$ in $  \mathfrak{\widetilde{G}} $ under the action of an inner automorphism $ Ad \xi $, $ \xi \in  \mathfrak{\widetilde{B}} $, $ \Phi ( \mathfrak{h})=(Ad \xi)\mathfrak{\widetilde{h}} $.

Let us introduce the element $ y=exp(\xi) $ and the subgroup $ \widetilde{G_{O}}=yG_{O}y^{-1} $. The analytic loop $ B(\star) $
$$
a \star b=\widetilde{\prod}(ab), a,b \in B
$$

where $ \widetilde{\prod}: G \longrightarrow B $ is the projection on $ B $ parallel to  $ \widetilde{G_{O}} $, is a Bol loop whose tangent Bol algebra is isomorphic to  $ \widetilde{\mathfrak{B}} $. In particular, $ B( \star) $ and   $ \widetilde{B}(\circ) $ are locally isomorphic. The operations on $ B $ are isotopic. Indeed let us introduce an analytic diffeomorphism $ \Omega: B \longrightarrow B $, $ a \longmapsto (y \times a) \setminus y $, then:\\
$$
\Omega^{-1}_{y}(\Omega_{y}(a) \star \Omega_{y}(b))=(L_{y})^{-1}(\Omega_{y}(a) \times L_{y}b)=Y^{-1}\times [((y\times a)\setminus y) \times (y \times b)]
$$
$$
=[y^{-1} \times (((y \times a) \setminus y) \times y^{-1})]\times (y^{2} \times b)
$$
$$
=[y^{-1} \times ((y^{-1} \times ((y^{2}\times a)\times y^{-1}))]\times(y^{2} \times b)
$$
$$
=[y^{-1} \times [y^{-1} \times ((y^{2}\times a)\times y^{-2})]\times(y^{2} \times b)
$$
$$
=[y^{-2} \times ((y^{2} \times a) \times y^{-2})]\times (y^{2} \times b)
$$
$$
=(a \setminus y^{2}) \times (y^{2} \times b)
$$
$$
=(a \setminus y) \times (y \times b)
$$
therefore 

$$ \Omega_{y}((a \setminus y) \times (y \times b))= \Omega_{y}(a) \star \Omega_{y}(b). $$
hence $B(\times)$ is isotopic to  $ \widetilde{B}(\circ) $ and by the 
diffeomorphism $\Omega$ they are isomorphic. 

Hence the theorem is proved.

{\bf Remarque:} The notion of isotopy of Bol algebra is not an equivalent 
relation in the manifold of Bol algebras. That is the reason why 
we say {\bf B- Isotopic Bol Algebra} instead of saying 
{\bf Isotopic Bol Algebra}

\section{ABOUT THE CLASSIFICATION OF BOL ALGEBRAS}

Let   $ \mathfrak{G}=\mathfrak{B} \dotplus \mathfrak{h} $ be Lie algebra in the
 involutive decomposition $  \mathfrak{h'} $ a subalgebra in  $ \mathfrak{G} $ 
such that  $ \mathfrak{G}=\mathfrak{B} \dotplus \mathfrak{h'} $
$$
a \cdot b=\prod [a,b]
$$
$$
(a,b,c)=[[a,b],c]
$$
where $ \prod:   \mathfrak{G} \longrightarrow \mathfrak{B}  $ is a projection 
on
  $ \mathfrak{B}  $ parallel to  $ \mathfrak{h'} $ and [,] commutator in Lie
 algebra  $ \mathfrak{G} $.\\
That is the raison why the classification of Bol algebra lead to the classification of subalgebra  $ \mathfrak{h'} $, in the enveloping Lie algebra  $ \mathfrak{G} $ (not necessarily canonical) for a Lie triple system  $ \mathfrak{M} $.\\
Below, we will examine the classification of Bol algebras with isomorphism accuracy and isotopic accuracy \cite{boumik}. The classification with isotopic accuracy is more crude, than the classification with isomorphism accuracy. However the notion of isotopy of Bol algebra is opening a new connection between non isomorphic Bol algebras.
\section{ISOCLINE BOL ALGEBRAS}

One can prove that any vectorial space  $ \mathfrak{B} $, equipped with 
operation
$$
\forall \xi, \eta, \zeta \in \mathfrak{B}\; \; \xi \cdot \eta=\alpha (\xi)\eta-\alpha(\eta)\xi, 
$$
$$
<\xi, \eta, \zeta>=\beta(\xi,\zeta)\eta-\beta(\eta, \zeta)\xi \; \; (I)
$$
where $ \alpha: \mathfrak{B} \longrightarrow \mathbb{R}$-linear form,  $ \beta: \mathfrak{B}\times \mathfrak{B} \longrightarrow \mathbb{R}$-bilinear symmetric 
form. Is a Bol algebra.

\underline{{\bf Definition 3.}} Bol algebra of view (I) is called isocline.

One can prove see \cite{sabmik1}, that any Bol algebra is called isocline if 
and only if
 it verifies the plane axiom that means any two-dimensional vectorial subspace 
(plane) is a subalgebra.

In particular any two-dimensional Bol algebra is isocline.

If $ dim \mathfrak{B}=3 $ and $ \alpha=0 $, then depending on the rank, and 
the 
signature of the form $ \beta $, we obtained 5 non trivial and non isomorphic 
types of Lie triple systems.

\section{BOL ALGEBRA WITH TRILINEAR OPERATION OF DEFINE BY $(e_2,e_3,e_2)=e_1 $ ; \;  $(e_2,e_3,e_3)=\pm e_2 $ }
 In what follows we consider Bol algebras of dimension 3, for their
construction see \cite{mik2, mike6}. In this paper we base our investigation 
of 
3-dimensional Bol algebras, on the examination of their canonical enveloping 
Lie algebras. As we already state it follows that the dimension of their 
canonical enveloping Lie algebras can not be more than 6. Below we limit 
ourselves to the classification of Bol algebras ( and their corresponding 
3-Webs), with canonical enveloping Lie algebras of dimension $ \leq 4$.

 Let us examine the case   $ dim \mathfrak{G}=4$, the structural constants of
 Lie algebra  $ \mathfrak{G}=<e_1,e_2,e_3,e_4>$, $\mathfrak{B}=<e_1,e_2,e_3>$ 
are defined as follows:
$$
[e_2,e_3]=e_4, \; \; [e_2,e_4]=- e_1
$$
\; \; \; \; \; \; \; \; \; \; \; \; \; \; \; \; \; \; \; \; \; \; \; \; \; \; (1)
$$
[e_3,e_4]=\mp e_2.
$$
In addition $ \mathfrak{G}=\mathfrak{B} \dotplus [\mathfrak{B},\mathfrak{B}]$,
$[\mathfrak{B}, \mathfrak{B}]=<e_4>$.

By introducing in consideration the 3-dimensional subspaces of subalgebras
$$
\mathfrak{h}_{x,y,z}=<e_4 +xe_1 +ye_2 +ze_3>,\; \; x,y,z, \in \mathbb{R}
$$
we obtain a collection of Bol algebras of view:
$$
e_2 \cdot e_3=-xe_1 -ye_2 -ze_3 ,
$$
\; \; \; \; \; \; \; \; \; \; \; \; \; \; \;\;\;\;\;\;\; \; \; \; \; \; \; \; (2)
$$
(e_2 ,e_3 ,e_2)=e_1,\;
$$
$$
 (e_2 ,e_3 ,e_3)=\pm e_2.
$$
\; Our main problem will be to give an isomorphical and isotopical
classification of Bol algebras of view (2).

\; For the complete examination of this case, we will split it in cases
Type $ V^-$ and Type $V^+$, corresponding to the upper and the lower signs 
of the formulas (1) and (2).

\; The group  $F$ of automorphisms of Lie triple system $\mathfrak{B}$
relatively to a fixed base $ e_1, e_2, e_3 $ from Type $V^-$ is defined as it
follows:

\begin{displaymath}
 F=\left\{A=\left(\begin{array}{ccc}
\pm b^2 & fb & d \\
0 & b & f \\ 
0 & 0 & \pm 1 \\
\end{array}\right); b \neq 0,  \right\}\; \; \;\;\;\;\;\;\;\;\;(3)
\end{displaymath}

\; The extension of an automorphism $A \in F$ to an automorphism of Lie
algebra $\mathfrak{G}$, transforming the subspace $\mathfrak{B}$ into itself
can be realized as it follows:

$$
Ae_4=A[e_2,e_3]=[Ae_2,Ae_3]= \pm be_4 . \; \; \;\;\;\;\;\;\; \; (4)
$$

In addition
$$
A(e_4 +xe_1 +ye_2 +ze_3)=\pm b \left(e_4 + \frac{xb^2 \pm yfb \pm zd}{b}e_1 \pm  \frac{yb+ zf}{b}e_2 +\frac{z}{b}e_3 \right),
$$
that is
$$
A(\mathfrak{h}_{x,y,z})=\mathfrak{h}_{x',y',z'},
$$
where
$$
x'= \frac{xb^2 \pm yfb \pm zd}{b},
$$
$$
y'=\pm  \frac{yb+ zf}{b} , \; \; \; \; \; \;\;\;\;\;\;\;\;\;(5)
$$
$$
z'=\frac{z}{b}.
$$

\begin{itemize}
\item If $ z \neq 0$, then by choosing $b, d$ and $f$ one can make 
 such that $(x,y,z)  \longrightarrow (x',y',z')  $: $x'=0,y'=0,z'=0$;
\item if $z=0$, hence $z'=0$ one can choose $y'= \pm y \geq 0$, and make
$x'=0$.
\end{itemize}
\; \; In this way, we obtain one family and one exceptional Bol algebra.

\underline{{\bf Proposition $ V^-$ .1}} Any Bol algebra of dimension 3, with the 
trilinear
operation of Type $V^-$ and the canonical enveloping Lie algebra of dimension 
4, is isomorphic to one of the following Bol algebras:
\begin{itemize}
\item $ V^- .1. \; \; e_2 \cdot e_3 =-e_3, \; (e_2,e_3,e_2)=e_1$, $(e_2,e_3,e_3)=e_2, $
\item $ V^- .2. \; \; e_2 \cdot e_3 =-ye_2 , \; (e_2,e_3,e_2)=e_1, (e_2,e_3,e_3)=e_2,
 \;  y \geq 0 $.
\end{itemize}
      
The distinguished Bol algebras are not isomorphic each to other.

Similarly one can establish the correctness of the following Proposition.

\underline{{\bf Proposition $ V^+$ .2}} Any Bol algebra of dimension 3, with the 
trilinear
operation of Type $V^+$ and the canonical enveloping Lie algebra of dimension 
4, is isomorphic to one of the following Bol algebras:
\begin{itemize}
\item $ V^+ .1.\; \; e_2 \cdot e_3 =-e_3, \; (e_2,e_3,e_2)=e_1$, $(e_2,e_3,e_3)=e_2 ,$
\item $ V^+ .2. \; \;e_2 \cdot e_3 =-ye_2, \; (e_2,e_3,e_2)=e_1, (e_2,e_3,e_3)=-e_2 
, \;y \geq 0. $ 
\end{itemize}
      
Also this distinguished Bol algebras are not isomorphic each to other.

Let us  pass to the isotopic classification of Bol algebras given in
Proposition $V^-$ .1.

We note that for every $ \xi=ue_1 +ve_2 +pe_3 $ from $ \mathfrak{B}$
$ u,v,p, \in \mathbb{R} $

\begin{displaymath}
ad ( \xi )= \left(\begin{array}{cccc}
0 & 0 & 0 & -v \\
0 & 0 & 0 & -p  \\
0 & 0 & 0 & 0 \\
0 & -p & v & 0 \\
\end{array}\right),
\end{displaymath}  

\begin{displaymath}
Ad ( \xi )= \left(\begin{array}{cccc}
1 & \frac{v(\cosh p -1)}{p} &- \frac{(\cosh p -1)v^2}{p^2} & -\frac{v \sinh p}{p} \\
0 & \cosh p & \frac{v(\cosh p -1)}{p} & -\sinh p  \\
0 & 0 & 1 & 0 \\
0 & -\sinh p & \frac{v \sinh p}{p} & \cosh p \\
\end{array}\right).
\end{displaymath}  
Let us find the the image of $\Phi(\mathfrak{h})$ under the action of
$ \Phi=Ad \xi$ on the one-dimensional subalgebra $\mathfrak{h}$ with a 
direction vector
$e_4  +ye_2,\; y \geq 0$; 

$$
\Phi (e_4  +ye_2)=(\cosh p -y\sinh p)e_4 +\frac{v}{p}[y(cosh p -1)-\sinh p]e_1 +(ycosh p -\sinh p)e_2, 
$$
 
$$
\Phi (e_4  +ye_2)=e_4 +\frac{v}{p}[\frac{y(cosh p -1)-\sinh p}{\cosh p -y\sinh p}]e_1 +\frac{ycosh p -\sinh p}{\cosh p -y\sinh p}e_2, 
$$
$$
x'=\frac{v}{p}[\frac{y(cosh p -1)-\sinh p}{\cosh p -y\sinh p}],
\; \; \; \; \; \; \; \; \; \; \; \; \; \; \; \; \; \; \; \; \; \; \;(5)
$$
$$
y'=+\frac{ycosh p -\sinh p}{\cosh p -y\sinh p}.
$$
\; \; By choosing $p \neq 0$ such that $y=\frac{\sinh p}{\cosh p -1}$,
we obtain $x'=0$. Let us note in addition that $coth p \neq 0$ that is the 
map is correctly defined. Applying to the obtain Bol algebra the 
automorphism of view 5, one can make $y'=1$.
\begin{itemize}
\item If $p =0$ then $x'=-v$\\
\;\;\;\;\;\;\;\;\; $y'=y$; \\
\item if $v \neq 0$ then one can make $ x'=1, y'=1$;\\
\item if $v=0$ then $x'=0$ and, one can make $y'=1$;\\
\item if $y=0$ then $x'=-\frac{v}{p}tanh p$,\\
\; \; \; \; \; \; \; \; \; \; \; \; \; \; \; \; \; \; \; \; 
\; \; \; \; \; \; \; \; (6)
 
$y'=-tanh p$,
\end{itemize}
\begin{itemize}
\item a) if $p=0$, then $x'=-v, v'=0$
\item b) if $p \neq 0$ then when applying the automorphism of view (6) and with
regards to $v$ we can make $x'=0; 1 , y'=1$.
\end{itemize}
\; \; The selected three cases are isotopic (in the sense of the definition 
of isotopy).

 We note the exceptional Bol algebra of Proposition $V^-$ under the 
action of isotopic transformation is not changing

\; \; Summarizing the conducted examination one can formulate the Theorem: 

\underline{{\bf Proposition $ V^-$ .3}} Any Bol algebra of dimension 3, with the 
trilinear
operation of Type $V^-$ and the canonical enveloping Lie algebra$ \mathfrak{G}$ of dimension 4, is isotopic to one of the following Bol algebras:
\begin{itemize}
\item $  e_2 \cdot e_3 =-e_3, \; (e_2,e_3,e_2)=e_1, (e_2,e_3,e_3)=e_2$;
\item $ e_2 \cdot e_3 =-e_2, \; (e_2,e_3,e_2)=e_1, \; (e_2,e_3,e_3)=e_2$;
\item $  e_2 \cdot e_3 =-e_1 -e_2, \; (e_2,e_3,e_2)=e_1, \; (e_2,e_3,e_3)=e_2$;
\item trivial bilinear operation, $(e_2,e_3,e_2)=e_1, \; (e_2,e_3,e_3)=e_2$. 
\end{itemize}
      
Analogically one can state the correctness of the Theorem.

\underline{{\bf Proposition$ V^+$ .4}} Any Bol algebra of dimension 3, with the 
trilinear
operation of Type $V^+$ and the canonical enveloping Lie algebra$ \mathfrak{G}$ of dimension 4, is isotopic to one of the following Bol algebras:
\begin{itemize}
\item $  e_2 \cdot e_3 =-e_3, \; (e_2,e_3,e_2)=e_1, (e_2,e_3,e_3)=-e_2$;
\item $ e_2 \cdot e_3 =e_2, \; (e_2,e_3,e_2)=e_1, \; (e_2,e_3,e_3)=-e_2$;
\item trivial bilinear operation, $(e_2,e_3,e_2)=e_1, \; (e_2,e_3,e_3)=-e_2$. 
\end{itemize}

\; \; Below, we reduce to description of 3-Webs corresponding to the isolated
Bol algebras of Type$ V^-$ and Type$ V^+$. 

The composition law $(\triangle)$,
corresponding to the Lie group $G$ of enveloping Lie algebra for Bol algebra is
defined as follows:

$$ \begin{bmatrix} x_1\\x_2\\x_3\\x_4\end{bmatrix} \triangle
   \begin{bmatrix} y_1\\y_2\\y_3\\y_4\end{bmatrix}
   = 
   \begin{bmatrix}
     x_{1}+y_{1}+ \frac{x_{4}y_{2} -y_{4}x_{2}}{2}\\x_{2}+y_{2}\cos x_{3}-y_{4}\sin x{3}\\x_{3}+y_{3}\\x_{4}-y_{2}\sin (x_3) +y_{4}\cos (x_3)\end{bmatrix}.
$$
 
 In case $ V^-.1$ the subgroup $ H=\exp \mathfrak{h} $, can be realized as the 
collection of elements

$$
H=\exp \mathfrak{h}=\{exp \alpha (e_4  +e_3)\}_{ \alpha \in \mathbb{R}}=\{0,0,\alpha, \alpha \}_{ \alpha \in \mathbb{R}}.
$$
The collection of elements
$$
B=\exp \mathfrak{B}=\left\{t+\frac{(v-\sin v)u^2}{2v^2},\frac{u}{v}\sin v,v,\frac{u}{v}(1-cos v)\right\}_{t,u,v \in \mathbb{R}}
$$
form a local section of left space coset $G \bmod H $.
$\exp: \mathfrak{G}\supset\mathfrak{B} \longrightarrow B \subset G$
and
$B =\exp \mathfrak{B}$
$$
\exp^{-1} \begin{bmatrix} x_1\\x_2\\x_3\\x_4\end{bmatrix}= \begin{bmatrix} x_1-\frac{(x_{2})^2 \sin^2 (x_3)}{2(x_3)^3}+\frac{(x_{2})^2 
\sin^3 (x_3)}{2(x_3)^4}\\\frac{x_2}{x_3}\sin x_3\\x_3\\x_4\end{bmatrix}
$$
$  x_1,x_2,x_3,x_4 \in \mathbb{R}$.

Any element $ (x_1,x_2, x_3, x_4) \in G $, in the neighborhood  $e$, can be
uniquely represented as follows:  
 
$$ \begin{pmatrix} x_1\\x_2\\x_3\\x_4\end{pmatrix} =
   \begin{bmatrix}x_1 +\frac{(x_3 -v)[x_2 +(x_3 -v)\sin v ]}{2}\\x_2 +(x_3 -v)\sin v\\ v\\x_4 -(x_3 -v)\cos v\end{bmatrix}
   \triangle 
   \begin{bmatrix} 0\\0\\x_3 -v\\ x_3 -v\end{bmatrix},
$$

where $v$ are any numbers defined from the relation
$$
\left[ x_4 -(x_3 -v)\cos v \right]\sin v =\left[x_2 +(x_3 -v)\cos v\right](\cos v -1).
$$

\; \; The composition law $( \star )$ corresponding to the local analytical Bol 
loop $B( \star )$ is defined as follows:
 
$$ \begin{pmatrix} t\\u\\v\end{pmatrix} \star
   \begin{pmatrix}t'\\u'\\v'\end{pmatrix}
   = \exp^{-1} \left(\prod_{B}\left(
    \begin{bmatrix}t\\u\\v\\0\end{bmatrix}\triangle \begin{bmatrix}
     t'\\u'\\v'\\0\end{bmatrix}\right) \right)
$$
$$ \begin{pmatrix} t\\u\\v\end{pmatrix} \star
   \begin{pmatrix}t'\\u'\\v'\end{pmatrix}
   = \exp^{-1} \left(\prod_{B}\left(
    \begin{bmatrix}t+t'\\u+u'\cos v\\v+v'\\u \sin (v)\end{bmatrix}\right) \right)
$$
$$
= \exp^{-1} \left( \begin{bmatrix}t+t'+\frac{(v+v'-T)\left[u+u'\cos v +(v+v'-T)\sin T\right]}{2}\\u+u'\cos v +(v+v'-T)\sin T\\T\\u\sin v -(v+v'-T)\cos T\end{bmatrix}\right)
$$.

$$
= \begin{bmatrix}F_1(t,t',u'u',v'v',T)\\\left[u+u'\cos v +(v+v'-T)sin T\right]\frac{\sin T}{T}\\T\end{bmatrix},
$$

where $T$ is defined from the relation:
$$
\left[ u\sin v -(v+v' -T)\cos T \right]\sin T =\left[u+u'\cos v +(v+v' -T)\cos T\right](\cos T -1).
$$

And $F_1(t,t',u'u',v'v',T)$ from the relation
\begin{multline}
F_1(t,t',u'u',v'v',T)=t+t'+\frac{(v+v'-T)(u+u'\cos v)}{2}+\frac{(v+v'-T)^2}{2}\sin T-\\
-\frac{\left[u+u'cos v + (v+v'-T)\sin T \right]^2}{2T^4}(T-\sin T)\sin^2 T.
\end{multline}

In  case $ V^-.2$ the subgroup $ H=\exp \mathfrak{h} $, can be realized as the 
collection of elements

$$
H=\exp \mathfrak{h}=\{exp \alpha (e_4 +ye_2 )\}_{ \alpha \in \mathbb{R}}=\{0,0, \alpha , \alpha \}_{ \alpha \in \mathbb{R}}.
$$
The collection of elements
$$
B= exp \mathfrak{B}=\left\{t+\frac{(v-\sin v)u^2}{2v^2},\frac{u}{v}\sin v,v,\frac{u}{v}(1-cos v)\right\}_{t,u.v \in \mathbb{R}}
$$
form a local section of left space coset $G \bmod H $.

Here $\exp^{-1}$  is defined as in the case above.

\; \; Any element $ (x_1,x_2, x_3, x_4) \in G $, in the neighborhood 
 $e$, can be uniquely represented as follows:  

$$ \begin{pmatrix} x_1\\x_2\\x_3\\x_4\end{pmatrix} =
   \begin{bmatrix}x_1 AB\frac{y-y\cos x_3 -\sin x_3}{ 2x_3}\\A\frac{\sin x_3}{x_3}\\ x_3\\A\frac{1-\cos x_3}{x_3}\end{bmatrix}
   \triangle 
   \begin{bmatrix} 0\\y\frac{x_{4}\sin x_3 -x_{2}(1-\cos x_3)}{y-y\cos x_3 +\sin x_3}\\0\\\frac{x_{4}\sin x_3 -x_{2}(1-\cos x_3)}{y-y\cos x_3 +\sin x_3}\end{bmatrix},
$$

where $A,B$ are any numbers defined from the relations
$$
A=x_3 \frac{x_2 (y\sin x_3 +\cos x_3)-x_4 (y \cos x_3 -\sin x_3)}{y-y\cos x_3 + \sin x_3},
$$
$$
B= \frac{x_4 \sin x_3 -x_2 (1- \cos x_3)}{y-y\cos x_3 + \sin x_3}.
$$ 

\; \; The composition law $( \star )$ corresponding to the local analytical Bol 
loop $B( \star )$ is defined as follows:
 
$$ \begin{pmatrix} t\\u\\v\end{pmatrix} \star
   \begin{pmatrix}t'\\u'\\v'\end{pmatrix}
   = \exp^{-1} \left(\prod_{B}\left(
    \begin{bmatrix}t\\u\\v\\0\end{bmatrix}\triangle \begin{bmatrix}
     t'\\u'\\v'\\0\end{bmatrix}\right) \right)
$$
$$ \begin{pmatrix} t\\u\\v\end{pmatrix} \star
   \begin{pmatrix}t'\\u'\\v'\end{pmatrix}
   = \exp^{-1} \left(\prod_{B}\left(
    \begin{bmatrix}t+t'\\u+u'\cos v\\v+v'\\-u\sin (v)\end{bmatrix}\right) \right)
$$

$$
= \begin{bmatrix}t+t'-C'\\A'\frac{\sin^2 (v+v')}{(v+v')^2}\\v+v'\end{bmatrix},
$$

where $A',C'$ are defined from the relations:
$$
A=(v+v') \frac{\left\{(u+u'\cos v)\left[ (y\sin (v+v') +\cos (v+v')\right]-u\sin v [y \cos (v+v') -\sin (v+v')]\right\}}{y-y\cos (v+v') + \sin (v+v')},
$$

$$
B= \frac{u\sin v \sin (v+v') -(u+u'\cos v) (1- \cos (v+v'))}{y-y\cos (v+v') + \sin (v+v')},
$$ 
$$
C'=A'B'\frac{(y-y\cos (v+v') -\sin (v+v'))}{v+v'}+(A')^2 \frac{\sin^4 (v+v')}{2(v+v')^4} \cdot \frac{-1+\cos (v+v')}{v+v'}.
$$

\; \; We pass to the description of Bol 3-Webs corresponding to the Type $V^+$.

\; \; In this case the composition law $(\triangle)$, corresponding to Lie group $G$, with 
enveloping Lie algebra of Bol algebra of Type $V^+$ is defined as:

$$ \begin{bmatrix} x_1\\x_2\\x_3\\x_4\end{bmatrix} \triangle
   \begin{bmatrix} y_1\\y_2\\y_3\\y_4\end{bmatrix}
   = 
   \begin{bmatrix}
     x_{1}+y_{1}+ \frac{x_{4}y_{2} -y_{4}x_{2}}{2}\\x_{2}+y_{2}\cosh x_{3}-y_{4}\sinh x{3}\\x_{3}+y_{3}\\x_{4}-y_{2}\sinh (x_3) +y_{4}\cosh (x_3)\end{bmatrix}.
$$
 
 In case $ V^+.1$ the subgroup $ H=exp \mathfrak{h} $, can be realized as the 
collection of elements

$$
H=\exp \mathfrak{h}=\{exp \alpha (e_4  +e_3)\}_{ \alpha \in \mathbb{R}}=\{0,0,\alpha, \alpha \}_{ \alpha \in \mathbb{R}}.
$$
The collection of elements
$$
B= \exp \mathfrak{B}=\left\{t+\frac{(v-\sinh v)u^2}{2v^2},\frac{u}{v}\sinh v,v,\frac{u}{v}(1-cosh v)\right\}_{t,u.v \in \mathbb{R}}
$$
form a local section of left space coset $G \bmod H $.
$\exp: \mathfrak{G}\supset\mathfrak{B} \longrightarrow B \subset G$
and
$B =\exp \mathfrak{B}$
$$
\exp^{-1}  \begin{bmatrix} x_1\\x_2\\x_3\\x_4\end{bmatrix}= \begin{bmatrix} x_1-\frac{(x_{2})^2 \sinh^2 (x_3)}{2(x_3)^3}+\frac{(x_{2})^2 \sinh^3 (x_3)}{2(x_3)^4}\\\frac{x_2}{x_3}\sinh x_3\\x_3\end{bmatrix}
$$
$  x_1,x_2,x_3,x_4 \in \mathbb{R}$.

Any element $ (x_1,x_2, x_3, x_4) \in G $, in the neighborhood  $e$, can be
uniquely represented as follows:  

$$ \begin{pmatrix} x_1\\x_2\\x_3\\x_4\end{pmatrix} =
   \begin{bmatrix}x_1 +\frac{(x_3 -v)[x_2 +(x_3 -v)\sinh v]}{2}\\x_2 +(x_3 -v)\sin v\\ v\\x_4 -(x_3 -v)\cos v\end{bmatrix}
   \triangle 
   \begin{bmatrix} 0\\0\\x_3 -v\\ x_3 -v\end{bmatrix}
$$

where $v$ are any numbers defined from the relation
$$
[ x_4 -(x_3 -v)\cosh v ]\sinh v =[x_2 +(x_3 -v)\cosh v](\cosh v -1).
$$

\; \; The composition law $( \star )$, corresponding to the local analytical Bol 
loop $B( \star )$, is defined as follows:
 
$$ \begin{pmatrix} t\\u\\v\end{pmatrix} \star
   \begin{pmatrix}t'\\u'\\v'\end{pmatrix}
   = \exp^{-1} \left(\prod_{B}\left(
    \begin{bmatrix}t\\u\\v\\0\end{bmatrix}\triangle \begin{bmatrix}
     t'\\u'\\v'\\0\end{bmatrix}\right) \right)
$$
$$ \begin{pmatrix} t\\u\\v\end{pmatrix} \star
   \begin{pmatrix}t'\\u'\\v'\end{pmatrix}
   = \exp^{-1} \left(\prod_{B}\left(
    \begin{bmatrix}t+t'\\u+u'\cosh v\\v+v'\\u \sinh (v)\end{bmatrix}\right) \right)
$$
$$
= \exp^{-1} \left( \begin{bmatrix}t+t'+\frac{(v+v'-P)\left[u+u'\cosh v +(v+v'-P)\sinh P\right]}{2}\\u+u'\cosh v +(v+v'-P)\sin P\\P\\u\sinh v -(v+v'-P)\cosh P\end{bmatrix}\right)
$$.

$$
= \begin{bmatrix}F_2(t,t',u'u',v'v',P)\\\left[u+u'\cosh v +(v+v'-P)sinh P\right]\frac{\sinh P}{P}\\P\end{bmatrix}
$$

where $P$ is defined from the relation:
$$
\left[ u\sinh v -(v+v' -P)\cosh P \right]\sinh P =\left[u+u'\cosh v +(v+v' -P)\cosh P\right](\cosh P -1)
$$

and $F_1(t,t',u'u',v'v',P)_1$ from the relation
\begin{multline}
F_1(t,t',u'u',v'v',P)=t+t'+\frac{(v+v'-P)(u+u'\cosh v)}{2}+\frac{(v+v'-P)^2}{2}\sinh P-\\
-\frac{\left[u+u'cosh v + (v+v'-P)\sinh P \right]^2}{2P^4}(P-\sinh P)\sinh^2 P.
\end{multline}

In case $ V^+.2$ the subgroup $ H=\exp \mathfrak{h} $, can be realized as the 
collection of elements

$$
H=\exp \mathfrak{h}=\{\exp \alpha (e_4 +ye_2 )\}_{ \alpha \in \mathbb{R}}=\{0,0, \alpha , \alpha \}_{ \alpha \in \mathbb{R}}.
$$
The collection of elements
$$
B= \exp \mathfrak{B}=\left\{t+\frac{(v-\sinh v)u^2}{2v^2},\frac{u}{v}\sinh v,v,\frac{u}{v}(1-cosh v)\right\}_{t,u.v \in \mathbb{R}}
$$
form a local section of left space coset $G mod H $.

Here $\exp^{-1}$  is defined as in the case above.

\; \; Any element $ (x_1,x_2, x_3, x_4) $ from $ G $, in the neighborhood 
 $e$, can be uniquely represented as follows:  

$$ \begin{pmatrix} x_1\\x_2\\x_3\\x_4\end{pmatrix} =
   \begin{bmatrix}x_1 DE\frac{y-y\cosh x_3 -\sinh x_3}{ 2x_3}\\D\frac{\sinh x_3}{x_3}\\ x_3\\D\frac{1-\cosh x_3}{x_3}\end{bmatrix}
   \triangle 
   \begin{bmatrix} 0\\y\frac{x_{4}\sin x_3 -x_{2}(1-\cosh x_3)}{y-y\cosh x_3 +\sinh x_3}\\0\\\frac{x_{4}\sinh x_3 -x_{2}(1-\cosh x_3)}{y-y\cosh x_3 +\sinh x_3}\end{bmatrix},
$$

where $D,E$ are any numbers defined from the relations
$$
D=x_3 \frac{[x_2 (y\sinh x_3 +\cosh x_3)-x_4 (y \cosh x_3 -\sinh x_3)]}{y-y\cosh x_3 + \sinh x_3},
$$
$$
E= \frac{x_4 \sinh x_3 -x_2 (1- \cosh x_3)}{y-y\cosh x_3 + \sinh x_3}.
$$ 

\; \; The composition law $( \star )$ corresponding to the local analytical Bol 
loop $B( \star )$ is defined as follows:
 
$$ \begin{pmatrix} t\\u\\v\end{pmatrix} \star
   \begin{pmatrix}t'\\u'\\v'\end{pmatrix}
   = \exp^{-1} \left(\prod_{B}\left(
    \begin{bmatrix}t\\u\\v\\0\end{bmatrix}\triangle \begin{bmatrix}
     t'\\u'\\v'\\0\end{bmatrix}\right) \right)
$$
$$ \begin{pmatrix} t\\u\\v\end{pmatrix} \star
   \begin{pmatrix}t'\\u'\\v'\end{pmatrix}
   = \exp^{-1} \left(\prod_{B}\left(
    \begin{bmatrix}t+t'\\u+u'\cosh v\\v+v'\\-u\sinh (v)\end{bmatrix}\right) \right)
$$

$$
= \begin{bmatrix}t+t'-F'\\D'\frac{\sinh^2 (v+v')}{(v+v')^2}\\v+v'\end{bmatrix}
$$

where $D',F'$ are defined from the relations:
$$
D'=(v+v') \frac{\Lambda}{y-y\cosh (v+v') + \sinh (v+v')},
$$

\begin{multline}
\Lambda=\{(u+u'\cosh v)\left[ (y\sinh (v+v') +\cosh (v+v')\right]-\\-u\sinh v \left[y \cosh (v+v') -\sinh (v+v')\right]\},
\end{multline}
$$
E'= \frac{u\sinh v \sinh (v+v') -(u+u'\cosh v) (1- \cosh (v+v'))}{y-y\cosh (v+v') + \sinh (v+v')},
$$ 
$$
F'=D'E'\frac{(y-y\cosh (v+v') -\sinh (v+v'))}{v+v'}+(A')^2 \frac{\sinh^4 (v+v')}{2(v+v')^4} \cdot \frac{-1+\cosh (v+v')}{v+v'}.
$$

\end{document}